\documentclass[english]{elsart}
\usepackage{babel}
\usepackage{amssymb}
\usepackage{epsfig,graphicx,times}                    
\begin{document}            
\begin{frontmatter}
\title{A 
  bijection between the $d$-dimensional simplices with distances
  in $\{1,2\}$ and the partitions of \vspace{-10mm}$d+1$}
\author{Christian Haase}
\thanks{The first author was partially supported by NSF-grant
  DMS-0200740.}
\ead{haase@math.duke.edu}
\address{Duke University, Department of Mathematics, Durham, NC
  27708-0320, USA}
\author{Sascha Kurz}
\ead{sascha.kurz@uni-bayreuth.de}
\address{University of Bayreuth, Department of Mathematics, D-95440 Bayreuth, Germany}
\end{frontmatter}
Integral point sets are sets of $n$ points in the Euclidean space $\mathbb{E}^d$ with
integral distances between vertices, see \cite{Har} for a survey. We examined such point sets for
$n=d+1$ and received the following table of numbers of nonisomorphic integral simplices 
by computer calculations. Here we call the largest occurring distance the diameter of the point set.
\begin{table}[h]        
  \begin{center}
    \begin{tabular}{|c|c|c|c|c|c|c|c|}
      \hline
      diameter & $d=3$ & $d=4$ & $d=5$ & $d=6$ & $d=7$ & $d=8$ & $d=9$ \\
      \hline
      1 & 1 & 1 & 1 & 1 & 1 & 1 & 1  \\
      2 & 4 & 6 & 10 & 14 & 21 & 29 & 41  \\
      3 & 16 & 56 & 197 & 656 & 2127 & 6548 & 19130  \\
      4 & 45 & 336 & 3133 & 31771 & 329859 & 3336597 & 32815796  \\
      \hline
    \end{tabular}\\[1mm]
    Table 1. Number of nonisomorphic integral simplices by diameter and dimension.
  \end{center}    
\end{table}
Due to the triangle inequality the (d+1)-element vertex set is
partitioned into subsets of vertices having pairwise distance $1$,
whereas vertices of different subsets are at distance $2$. To prove
the proposed bijection, we have to provide a simplex for a given
partition $(n_1,\dots,n_r)$ of $d+1$. We would like to mention that
the bijection holds more generally for simplices with distances in
$\{1,\lambda\}$ for $\lambda\ge 2$. At first we give the following
explicit construction.

\textbf{Construction.} Place regular $(n_i-1)$-simplices with edge
length $1$ with their barycenters at the origin into mutually
orthogonal spaces. Then shift the $i^{{th}}$ simplex into a new
coordinate direction by the amount of 
$\sqrt{\frac{\lambda^2}{2}-\frac{n_i-1}{2n_i}}$.\\[-6mm]

For another proof we need the following criterion.\\[-6mm]

\textbf{Theorem (Menger \cite{Menger}).}
If $M$ is a set of $d+1$ points with distance matrix $D=(d_{i,j})$
and $A=(d^2_{i,j})$, then $M$ is realizable in the Euclidean
$d$-dimensional space, iff $(-1)^{d+1}\det(\overline{A})\ge 0$ and
each subset of $M$ is realizable in the $(d-1)$-dimensional space,
where 
$\overline{A}:=\left(
  \begin{array}{cc}
    0 & (1,\dots,1)^T \\
    (1,\dots,1) & A \\
  \end{array}\right).$
To apply this theorem we provide, for distance matrices derived from a
partition of $d+1$ and with the nonzero values being in
$\{1,\lambda\}$, the following lemma.

\textbf{Lemma.}
$\quad (-1)^{d+1}(\lambda^2\, \det(\overline{A}) + \det(A)) > 0 \quad
\mbox{and} \quad (-1)^{d+1} \det(\overline{A}) > 0.$

We leave the proof to the reader, because it can be easily but
lengthly done by induction on $d$. It should be remarked that a
formula for the relevant determinants was also stated in
\cite{small_diameter}, but with no details of the computation.\\[2mm] 
As a last remark we would like to mention that using [5] one can
generalize the stated bijection. For given $d$ only
$\lambda\ge\sigma(d-1,d+1)$ is needed, where
$$
\sigma(d,d+2) = \sqrt{\frac{9d-10+\sqrt{33d^2-52d+20}}{4d-4}} \ge
\frac{1}{2} \sqrt{9+\sqrt{33}} \approx 1.91993.
$$ 
\vspace{-10mm}
\bibliography{Bijection_Simplices_Partitions_2}

\begin{thebibliography}{1}

\bibitem{Partition}
G.~Andrews.
\newblock {\em The theory of partitions}, volume~2 of {\em Encyclopedia of
  Mathematics and its Applications}.
\newblock Cambridge University Press, 1984.

\bibitem{0548.51014}
A.~Blokhuis.
\newblock {\em {Few-distance sets.}}
\newblock {CWI Tracts, 7. Centrum voor Wiskunde en Informatica. Amsterdam:
  Mathematisch Centrum. IV, 70 p. Dfl. 10.80 }, 1984.

\bibitem{Har}
H.~Harborth.
\newblock {Integral distances in point sets.}
\newblock In {\em {Butzer, P. L. (ed.) et al., Karl der Grosse und sein
  Nachwirken. 1200 Jahre Kultur und Wissenschaft in Europa. Band 2:
  Mathematisches Wissen. Turnhout: Brepols. 213-224 }}. 1998.

\bibitem{small_diameter}
H.~Harborth and L.~Piepmeyer.
\newblock {Points sets with small integral distances.}
\newblock In {\em {Applied geometry and discrete mathematics, Festschr. 65th
  Birthday Victor Klee, DIMACS, Ser. Discret. Math. Theor. Comput. Sci. 4,
  319-324 }}. 1991.

\bibitem{0789.05035}
H.~Harborth and L.~Piepmeyer.
\newblock {Two-distance sets and the golden ratio.}
\newblock In {\em {Bergum, G. E. (ed.) et al., Applications of Fibonacci
  numbers. Volume 5: Proceedings of the fifth international conference on
  Fibonacci numbers and their applications, University of St. Andrews,
  Scotland, July 20-24, 1992. Dordrecht: Kluwer Academic Publishers. 279-288
  }}. 1993.

\bibitem{Menger}
K.~Menger.
\newblock {Untersuchungen \"uber allgemeine Metrik.}
\newblock {\em Math. Ann.}, 100:75--163, 1928.

\end{thebibliography}
\bibliographystyle{abbrv}
\nocite{0548.51014,0789.05035,Partition}
\end{document}